\documentclass[11pt]{amsart}

\usepackage{amssymb,eucal}
\usepackage{amscd}
\usepackage[all,cmtip]{xy}
\usepackage{rotating}
\usepackage{hyperref,mathrsfs}
\usepackage[osf]{mathpazo}
\usepackage{tikz}
\usepackage{rotating}

\usepackage{amsmath,amssymb}

\newtheorem{theorem}{Theorem }[section]

\newtheorem{lemma}[theorem]{Lemma}

\newtheorem{remark}[theorem]{Remark}
\newtheorem{corollary}[theorem]{Corollary}
\newtheorem{proposition}[theorem]{Proposition}

\newcommand{\proj}{\mathrm{proj}}
\def\I{\mathbf{I}}
\def\nI{\mathop{\not\mathrm{I}}}

\def\PG{\mathbf{PG}}
\def\PGL{\mathbf{PGL}}
\def\PSL{\mathbf{PSL}}

\def\eop{\hspace*{\fill}$\Box$}
\newcommand{\Aut}{\mathrm{Aut}}
\newcommand{\id}{\mathbf{1}}

\newcommand{\cL}{\mathscr{L}}

\newcommand{\mP}{\mathscr{P}}
\newcommand{\cQ}{\mathscr{Q}}

\newcommand{\mS}{\mathscr{S}}

\newcommand{\cP}{\mathscr{P}}
\newcommand{\hW}{\mathbf{W}}

\newcommand{\hT}{\mathbf{T}}
\newcommand{\mB}{\mathscr{B}}

\title[Criterion concerning Singer groups]{A criterion concerning Singer groups of generalized quadrangles, and construction of uniform lattices in
{\large $\widetilde{\mathbf{C}_2}$}-buildings}


\author{Stefaan De Winter and Koen Thas}

\email{sgdwinte@me.com} 
\address{Michigan Technological University, Department of Mathematics, Fisher Hall, 1400 Townsend Drive, Houghton MI49931, USA}
\urladdr{http://web.me.com/sgdwinte/Homepage$\_$of$\_$Stefaan$\_$De$\_$Winter}

\email{kthas@cage.UGent.be}
\address{{Ghent University},
{Department of Mathematics},
{Krijgslaan 281, S25, B-9000 Ghent, Belgium}}
\urladdr{http://cage.ugent.be/$\sim$kthas}

\thanks{The first author is a  Postdoctoral Fellow
of the Fund for Scientific Research --- Flanders (Belgium). The second author was partially supported by the Fund for Scientific Research --- Flanders (Belgium)
when the present research was conducted.}
\date{}

\begin{document}
\maketitle

\begin{abstract}
We describe a simple criterion to construct Singer groups of Payne-derived generalized quadrangles, yielding, as a corollary, a classification of Singer groups of the classical
Payne-derived quadrangles in any characteristic. This generalizes recent constructions of Singer groups of these quadrangles that were presented in \cite{BambergGiudici}. 
In the linear case, and several other cases, our classification is complete. Contrary to what seemed to be a common belief, 
we show that for the classical Payne-derived quadrangles,
the number of different Singer groups is extremely large, and even
bounded below by an exponential function of the order of the ground field.  \\
Our results have direct applications to the theory of $\widetilde{\mathbf{C}_2}$-buildings, which are explained at the end of the paper.
\end{abstract}

\setcounter{tocdepth}{1}
\tableofcontents

\bigskip
\section{Introduction}

The study of projective planes admitting a so-called ``Singer group'', that is, a group acting sharply transitively on its point set, has been a central topic in Finite Geometry and difference set theory since Singer's seminal paper \cite{Singer} on the existence of such groups for the geometries $\mathbf{PG}(2,q)$ (that is, projective planes defined over the finite field $\mathbb{F}_q$). It is now of course natural to wonder not only which planes apart from $\mathbf{PG}(2,q)$, if any, admit Singer groups, but also whether other (building-like) geometries can admit a Singer group. It was D. Ghinelli \cite{Ghinelli} who initiated a theory of Singer groups for generalized quadrangles (GQs), i.e., (spherical) buildings of type $\mathbf{B}_2$. Over the last five years the authors of the present paper, together with E. E. Shult, have contributed to the theory of Singer groups for GQs in a series of papers \cite{SDWKT1, SDWKT3,  SDWKT2, SDWEESKT1}. 
\medskip

In \cite{SDWEESKT1} the known GQs admitting a Singer group  were classified in a combinatorial fashion. Through a different and less elementary approach from the combinatorial point of view, and using the recent classification of regular subgroups of almost simple
primitive  groups \cite{RegPrimGroups}, J. Bamberg and M. Giudici \cite{BambergGiudici} also obtained the classification of classical GQs admitting a Singer group (we credit them for having pointed out a small error in \cite{SDWEESKT1} which led to one small example being overlooked in \cite{SDWEESKT1}). By using \cite{RegPrimGroups} they also obtained the classification of the possible Singer groups that can act on a classical GQ. In their paper they provide examples  of groups that act as a Singer group on the so-called ``Payne-derived GQ'' (cf. the next section) of the classical symplectic GQ $\hW(q)$.  Some of these latter examples seemed to contradict the main theorem of \cite{SDWKT2}. 
\medskip

Initially, the aim of this paper was to correct the statement of the main theorem of \cite{SDWKT2}; in the last line of the proof of that theorem an inaccuracy occurs. (However, the main result, namely that a GQ admitting an odd order Singer group of order $p^m$, $p$  prime and $n \in \mathbb{N}$, with a center of order at least $\sqrt[3]{p^m}$, has to be Payne-derived, remains valid.)  When writing up the correction, a transparent criterion was noticed by the authors which opened up the possibility to classify Singer groups of general Payne-derived GQs. For the classical case | that is, when one derives from the symplectic $\hW(q)$-quadrangle where $q$ is any prime power | the criterion becomes particularly handy.

\medskip 
The present goal is to provide a classification of the groups that can act as a Singer group on the Payne-derived GQ of $\hW(q)$, as a corollary of this criterion.
All constructions of \cite{BambergGiudici} follow easily without performing any calculation. (It should be mentioned that the constructions of the present paper were obtained before
the current version of \cite{BambergGiudici} was available.)

The classification we propose is complete in many cases, the most important one being the linear case (i.e., when the group is induced by a linear group of the ambient projective space). In fact, contrary to what seemed to be a common belief, 
the number of different Singer groups for the classical Payne-derived quadrangles
is {\em extremely} large, and even
bounded below by an exponential function of the order of the ground field. On the way observing this, we also count (in a precise way) the number of linear sharply transitive groups of the classical affine plane (defined over a finite field $\mathbb{F}_q$). In fact, all these results also hold for general finite semifield planes (that is, planes which are both translation and dual translation planes). 
To our knowledge, also these observations (and techniques) are new.\\

In a second part of the paper, we analyse Singer groups of the known nonclassical Payne-derived quadrangles (which all occur in even characteristic). 

Finally, using recent constructions of Jan Essert, we apply our results to the construction theory of uniform lattices in buildings of type $\widetilde{\mathbf{C}_2}$.\\
\medskip

{\em Acknowledgement}.\quad
We are grateful to J. Bamberg and M. Giudici for kindly pointing out the error in \cite{SDWKT2} and sharing several versions of  their preprint \cite{BambergGiudici} with us.
We also want to thank Jan Essert and Linus Kramer for several helpful suggestions.\\

\section{Payne-derivation, and some notation}

Throughout this note we will use the same terminology and notation as introduced in \cite{SDWKT2}. (For definitions of the concepts related to GQs not defined here we also refer the reader to  \cite{PTsec}, or again \cite{SDWKT2}.)  

A {\em generalized quadrangle} (GQ) $\Gamma = (\mP_{\Gamma},\mB_{\Gamma},\I)$ is a rank $2$ geometry  of which the incidence graph is bipartite of girth $8$ and diameter $4$. 
We only consider {\em finite} GQs, meaning that $\vert \mP_{\Gamma} \vert$ (and then $\vert \mB_{\Gamma}\vert$) is finite.
(In our notation, $\mP_{\Gamma}$ is the point set and $\mB_{\Gamma}$ the line set.) When each line is incident with at least
$3$ points and each point is incident with at least $3$ lines, that is, when the geometry is ``thick'', 
it can easily be shown that there are constants $s$ and $t$ such that there are $s + 1$ points on each line and $t + 1$ lines on each point. We call $(s,t)$ the {\em order} of the quadrangle. If points $x$ and $y$ are collinear points, we write $x \sim y$ (and allow $x \sim x$), and $x^{\perp} = \{ z \in \mP_{\Gamma} \vert z \sim x \}$. For any $A \subseteq \mP_{\Gamma}$,

\begin{equation}
A^{\perp} =  \bigcap_{a \in A}a^{\perp} \ \ \mbox{and}\ \ A^{\perp\perp} = {(A^{\perp})}^{\perp}.
\end{equation}

If $x \not\sim y$, then $\vert \{x,y\}^{\perp}\vert = t + 1$, so that $\vert \{x,y\}^{\perp\perp}\vert \leq t + 1$. If equality holds for all $y \ne x$, we say that $x$ is {\em regular}.

We now quickly recall Payne's construction of GQs from  \cite{Payne71}. Let $\Gamma$ be a GQ of order $s$ (meaning of order $(s,s)$) with a regular point $x$. Define a rank $2$ geometry $\cP = \cP(\Gamma,x)$ as follows.

\begin{itemize}
\item The points of $\cP$ are the points of $\Gamma$ not collinear with $x$.
\item The lines of $\cP$ are the lines of $\Gamma$ not through $x$ together with all sets $\{x,r\}^{\bot\bot}\setminus\{x\}$, where $x\not\sim r$.
\item
Incidence is the natural one.
\end{itemize}

We have the following theorem.

\begin{theorem}[S. E. Payne \cite{Payne71}]
The point-line geometry $\cP$ is a generalized quadrangle of order $(s-1,s+1)$.
\end{theorem}

The GQ $\cP$ is the so-called {\em Payne-derived GQ} or {\em Payne-derivative} of $\Gamma$ (with respect to $x$), and sometimes we say that $\Gamma$ is ``Payne-integrated'' (with respect to $x$).
In the case where $\Gamma \cong \hW(q)$, every point of $\Gamma$ is regular \cite[Chapter 3]{PTsec}, and all Payne-derivatives are isomorphic, independent of the regular point chosen. We will denote the Payne-derivative of $\hW(q)$ therefore simply by $\cP(q)$.

The correct statement of the Main Theorem of \cite{SDWKT2} should read:

\begin{theorem}
\label{Main}
Let $\cQ$ be a GQ of order $(s,t)$ admitting a Singer group $G$, where $G$ is a $p$-group and $p$ is odd. Suppose $\vert Z(G)\vert \geq \sqrt[3]{\vert G\vert}$.
Then the following properties hold.
\begin{itemize}
\item[{\rm (1)}]
We have $t = s + 2$, and there is a GQ $\cQ'$ of order $s + 1$ with a regular point $x$, such that $\cQ$ is the Payne-derivative of $\cQ'$ with respect to $x$.
\item[{\rm (2)}]
We have $\vert Z(G)\vert = \sqrt[3]{\vert G\vert}$, that is,  $\vert Z(G)\vert = s + 1$.
\end{itemize}
\end{theorem}

The original version of the theorem as it appeared in \cite{SDWKT2} had the extra claim in (1) that ``The GQ $\cQ'$ is an STGQ of type $({\cQ'}^{x},K)$, with $K$ isomorphic to $G$.'' This is however not true, as the claim at the end of the proof that $G$ induces a full (= transitive) group of elations of $\cQ'$ is false.

The latter fact will become very clear in the present paper.\\

\section{A simple criterion}

In this section we describe a simple criterion to construct Singer groups for Payne-derived GQs in rather general circomstances.  

\subsection{Criterion}

Consider the following general situation.
Suppose $\mS$ is a thick GQ of order $s$, and let $x$ be a regular point of $\mS$.  Suppose $G$ is a Singer group of $\mP(\mS,x)$ with the following properties: 

\begin{itemize}
\item
it is induced by some automorphism group of $\mS$ | that is, by some subgroup $\overline{G}$ of $\mathrm{Aut}(\mS)_x$; 
\item
$\overline{G}$ contains a group $\mathbb{S}$ of order $s$ consisting of symmetries with center $x$.   
\end{itemize}

We mention the following lemma without further notice.

\begin{lemma}
$\mathbb{S} \unlhd \Aut(\mS)_x$. \eop \\
\end{lemma}

Construct the standard affine plane $\Pi(x)$ of order $s$ from the regular point $x$ \cite[1.3.1]{PTsec}; its points are the sets $\{x,z\}^{\perp}$ with $z \not\sim x$,
its lines are the elements of $x^{\perp} \setminus \{x\}$. Note that elements of $\mathrm{Aut}(\mS)_x$ induce elements of $\mathrm{Aut}(\Pi(x))$ via the following map

\begin{equation}
\xi: \mathrm{Aut}(\mS)_x \mapsto \mathrm{Aut}(\Pi(x)): g \mapsto g\mathbb{S}.
\end{equation}

Then $\overline{G}/\mathbb{S}$ is an automorphism group of $\Pi(x)$ which acts sharply transitively on its points. 

Vice versa, let $K$ be a subgroup of $\mathrm{Aut}(\Pi(x))$ which acts sharply transitively on its points, and suppose it is induced by some automorphism group $K^\#$ of $\mS$
(which then has to fix $x$).
Then $\overline{K} = \langle \mathbb{S},K^\#\rangle$, the group obtained by adjoining $\mathbb{S}$ to $K^\#$, induces a Singer group of $\mP(\mS,x)$ 
(note that $\vert \overline{K} \vert = s^3$, and that the kernel of the action of $\overline{K}$ on the points of $\Pi(x)$ is $\mathbb{S}$).\\

In particular, if $K$ is a translation group of $\Pi(x)$, then $\overline{K}$ is a (complete) elation group for $\mS^x$. (Vice versa, an elation group of $\mS^x$
containing $\mathbb{S}$ induces a translation group of $\Pi(x)$.) If $K$ is as such, $s$ is of course a prime power.

\subsection{Applications}

In general, automorphisms of Payne-derived GQs are not induced by automorphisms of the ambient GQ | cf. the discussion in \cite{SDWKT3}. 

For the classical case, we have the following satisfying answer:

\begin{theorem}[T. Grundh\"{o}fer, M. Joswig and M. Stroppel \cite{Strop}]
\label{Stropp}
If $q \geq 5$, any automorphism of $\mP(q)$ is induced by an automorphism of $\hW(q)$ fixing $x$.
\end{theorem}

All known (counter) examples of automorphisms not coming from the ambient GQ exist in even characteristic (besides some small sporadic examples in odd characteristic), see \cite{SDWKT3}.
The most general result available is the following:

\begin{theorem}[S. De Winter and K. Thas \cite{SDWKT3}]
\label{ind}
Let $\mP(\mS,x)$ be Payne-derived from the thick GQ $\mS$ of order $s$, with $s$ odd and $s \geq 5$. If $x$ is a center of symmetry, then any automorphism
of $\mP(\mS,x)$ is induced by an element of $\mathrm{Aut}(\mS)_x$. 
\end{theorem}

Let $\mS$ be a GQ satisfying the conditions of the previous theorem, and suppose $K$ is a Singer group of $\mP(\mS,x)$.
Then by Theorem \ref{ind}, and the criterion of the previous paragraph, $K$ is induced by an automorphism group $\overline{K}$ of $\mS$ of order $s^3$ which induces
a point-regular automorphism group of $\Pi(x)$ if it contains $\mathbb{S}$. We will encounter various general situations in which the latter assumption is 
naturally satisfied.\\


\section{Singer groups of $\cP(q)$}

We now consider  some concrete cases.

\subsection{Odd characteristic}

The only known thick GQ of order $s$ with $s$ odd having a regular point is the classical $\hW(q)$-GQ, with $q = s$. 
Let $\mS \cong \hW(q)$, $q$ odd, and let $x$ be any point. Also, we suppose that $q \geq 5$.\\

We note the following facts:
\begin{itemize}
\item
by the previous section, $\mathrm{Aut}(\mP(\hW(q),x)) \cong \mathrm{Aut}(\hW(q))_x$;
\item
$\Pi(x)$ is a Desarguesian affine plane;
\item
$\mathrm{Aut}(\Pi(x)) \cong \mathrm{Aut}((\hW(q))_x/\mathbb{S}$.
\end{itemize}

Let $\mathbf{S}(\Pi(x))$ denote  the set of sharply transitive groups of $\Pi(x)$, and $\mathbf{S}^{\mathbb{S}}(\mP(q))$ the set of Singer groups of $\mP(q)$ containing $\mathbb{S}$.
By the previous section and aforementioned properties, the next theorem stands.

\begin{theorem}
There is a natural bijection\\

\begin{equation}
\label{lab}  
\chi:  \mathbf{S}(\Pi(x))  \mapsto  \mathbf{S}^{\mathbb{S}}(\mP(q)) : K \mapsto \chi(K),  
\end{equation}

\medskip
\noindent
which sends $K$ to $\overline{K} \leq \mathrm{Aut}(\hW(q))_x$, and identifies $\overline{K}$ with the automorphism group $\chi(K)$ induced by $\overline{K}$ on 
$\mP(q)$\eop .\\
\end{theorem}

\begin{remark}
{\rm 
\begin{itemize}
\item[(i)]
Clearly isomorphic elements of $\mathbf{S}^{\mathbb{S}}(\mP(q))$ are mapped by $\chi^{-1}$ to isomorphic elements of $\mathbf{S}(\Pi(x))$. 
\item[(ii)]
If the conditions of Theorem \ref{Main} are satisfied, then  we still have a natural injection $\mathbf{S}^{\mathbb{S}}(\mP(\cQ',x)) \hookrightarrow  \mathbf{S}(\Pi(x))$, but since 
we do not know the precise relation between $\mathrm{Aut}(\Pi(x))$ and $\mathrm{Aut}(\cQ')_x$, we cannot claim a bijection as in (\ref{lab}).
\end{itemize}
}
\end{remark}

\medskip
By (\ref{lab}), the groups of Theorem \ref{Main} are completely classified when $\cQ' \cong \hW(q)$.\\

\medskip
\begin{theorem}
\label{thmoddcent}
Let $q$ be odd, $q \geq 5$.
If $T$ is a Singer group of $\mP(q)$ contained in $\mathbf{PGL}_4(q)_x$ (as a subgroup of $\mathrm{Aut}(\hW(q))_x$), then $\mathbb{S} \leq T$. 
\end{theorem}

{\em Proof}.\quad
We know that $T$, as a subgroup of $\mathrm{Aut}(\hW(q))_x$, fixes some line $L \I x$. It is clear that $T$ acts transitively on $L \setminus \{x\}$, so that for any
point $y \in L \setminus \{x\}$, $\vert T_y\vert = q^2$. Moreover, as $T_y \leq \mathbf{PGL}_4(q)_x$ and $T$ is a $p$-group with $q$ a power of the odd prime $p$, 
$T_y = K$ fixes $L$ pointwise (and hence is independent of the choice of $y$). No element of $K^{\times}$ can fix a line not meeting $L$ since $K \leq T$, so 
$K$ is a subgroup of the translation group of the TGQ $\hW(q)^L$. Now note that $K_{[x]}$, the subgroup of $K$ fixing $x$ linewise, has size $rq$ for some 
nonzero natural number $r$. Also, $K_{[x]}$ has the property that any of its elements either is a symmetry with center $x$, or an $(x,L,u)$-elation, with 
$u \I L \I x \ne u$ (such an $(x,L,u)$-elation fixes each of $x, L, u$ elementwise). (For a proof, use Benson's theorem, cf. \cite{PTsec}, \S 1.9. It is important to note that no nontrivial $(x,L,u)$-elation can fix another point $v$ on $L$
linewise besides $x$ and $u$ | this is because $L$ is antiregular.)
Let $a$ be the number of symmetries with center $x$ in $K_{[x]}$, including 
the trivial one. Since $T$ fixes $x$ and acts transitively on $L \setminus \{x\}$, the fact that $K_{[x]}$ is a normal subgroup of $T$ implies that the number 
of nontrivial  $(x,L,u)$-elations in $K_{[x]}$ is independent of the choice of $u$, and hence the cardinality of $\cup_{u \I L \I x \ne u}H(x,L,u)^{\times}$ is a multiple of $q$, where $H(x,L,u)$ is the group of $(x,L,u)$-elations. So 
\begin{equation}
rq - a \equiv 0\mod{q}.
\end{equation}
It follows that $a = q$, that is, $\mathbb{S} \leq T$.
\eop \\

\begin{remark}
{\rm Note that if $T$ is a Singer group of $\mP(q)$ which is contained in $\mathbf{PGL}_4(q)_x$ (as a subgroup of $\mathrm{Aut}(\hW(q))_x$), then $T$
automatically contains $\mathbb{S}$, so such a  group always comes, via $\chi$, from some element $\chi^{-1}(T)$ of   $\mathbf{S}(\Pi(x))$. Clearly, $\chi^{-1}(T)$ is a subgroup 
of $\mathbf{PGL}_3(q)_{[\infty]}$, where $[\infty]$ is the line at infinity of the projective completion of $\Pi(x)$.}\\
\end{remark}

\subsection*{Example: Translation group}

Let $K$ be the translation group of $\Pi(x)$ | it is elementary abelian, and fixes every point at infinity.
Then $\chi(K)$ is the Singer group of $\mP(q)$ which is induced by the standard elation group of $\hW(q)$ w.r.t. $x$.\\

\subsection*{Example: Abelian case}

More generally, suppose $K$ is abelian. Then by \cite{DemOst}, either $K$ is the translation group, and we can use the previous paragraph, or $K$
has precisely two orbits at infinity, namely one of size $1$ (a point, say $z$) and one of size $q$. Such a group is called a group ``of type (b)''.
The next theorem completely classifies these groups in {\em any} characteristic (still in the classical case), taken the discussion after the theorem into account.

\begin{theorem}[Classification of groups of type (b)]
Let $K$ be an abelian group which acts transitively on the points of the affine Desarguesian plane $\Pi(x)$ over $\mathbb{F}_q$, $q$ any prime power. Let $K$ be of type (b), fixing 
the point $z$ at infinity.
Extend $K$ by adjoining the translation group $T$ to $K$ | one obtains
a group $H$ of size $q^3$ which is generated by $T$ and the translation group $T^D$ w.r.t. $z$ in the dual setting. So $H$ is the classical semifield group (the Heisenberg group of dimension $3$ over $\mathbb{F}_q$), and
it is isomorphic to

\begin{equation} \{\left(
 \begin{array}{ccc}
 1 & \alpha & c\\
 0 & \mathbb{I} & \beta^T\\
 0 & 0 & 1\\
 \end{array}
 \right)\vert  \alpha, \beta \in \mathbb{F}_{q}, c \in \mathbb{F}_q \},                                         \end{equation}

\noindent
with standard matrix multiplication, and where $\mathbb{I}$ is the identity $2 \times 2$-matrix.
It follows that $K$ is a linear group.
\end{theorem}

{\em Proof}.\quad
It is clear that $K$ normalizes $T$, so $H$ is a group of order $q^3$ and $K$ and $T$ are abelian subgroups of order $q^2$, meeting in a subgroup $Z$ of size $q$
which fixes the line at infinity $[\infty]$ pointwise and $z$ linewise (the considerations being made in the plane). We emphasize that the fact  that $K$ acts sharply transitively on the points is crucial here.
It is immediate that $Z$ is the center of $H$ (since $K$ and $T$ are both abelian). 
(The inequality $K \cap T \leq Z$ is obvious. Conversely, if $z \in Z \setminus (K \cap T)$, then $\langle z,T \rangle$ is an abelian group which acts 
transitively but not sharply transitively on the points of $\Pi(x)$, contradiction.)
If $a \in T$ and $b \in K$ (not both simultaneously contained in $T$ nor $K$) are elements in $H$ for which
$[a,b] = \mathbf{1}$, then either $a \in Z$ or $b \in Z$.  
It follows that $C_K(T) = Z$ and $C_T(K) = Z$ (while $C_H(T) = T$ and $C_H(K) = K$). 

Let $k \in K$, $t \in T$, and suppose $[k,t] \not\in Z$. As $[k,t] \in T$, $[k,t]$ fixes some point $u \I [\infty]$ linewise,
      with $u \ne z$. Let $V \I u$, $V \ne [\infty]$, and put $U = V^{k^{-1}}$. Then 
      \begin{equation}
      U^{k^t} = {(V^{k^{-1}})}^{k^t} = V^{[k,t]} = V^{[t,k]} = V. 
      \end{equation}
      
      Write $t = t't''$, with
      $t' \in Z$ and $t'' \in T_{[u']}$, where $T_{[u']}$ is the subgroup of $T$ fixing $u'$ linewise, $u'$ being $u^{k^{-1}}$. (Note that
      $T = ZT_{[u']}$.) Then 
      \begin{equation}
      U^{k^t} = U^{k^{t't''}} = U^{k^{t''}} = U^{t''^{-1}kt''} = U^{kt''} = V^{t''}. 
      \end{equation}
      As the latter must be $V$, it follows
      that $t''$ fixes $V$, which is only possible if $t'' = \mathbf{1}$, that is $t \in Z$, contradiction. 
We have obtained that
\begin{equation}
[H,H] \leq Z = T \cap K,
\end{equation}
and so $H$ is of class $2$. Now take any $k \in K \setminus Z$; then
\begin{equation}
\vert \{ [k,t] \vert t \in T \} \vert = \vert T/C_T(k) \vert = \vert Z\vert 
\end{equation}
so that $[H,H] = Z$. 

It follows that  $H/Z$ is an abelian group, so that $H/N$
with $N$ the kernel of the action of $H$  on the lines incident with $z$, is a group of size $q$, and $N$ is a group of size $q^2$ which fixes $z$ linewise.
Whence $N = T^D$, and $H = \langle T,T^D\rangle$.  The theorem now follows from, e.g., \cite{Hira}.
\eop \\

It is well-known \cite{Hira} that $H$ contains $q + 1$ maximal (elementary) abelian groups of size of $q^2$ | including $T$, $T^D$ and $K$ | {\em now taken that $q$ is odd}. 
Let $M$ be the set of $q + 1$ maximal abelian groups (of order $q^2$) inside $H = \langle T,T^D\rangle$; note that
      elements of $M$ intersect two by two in $Z$, and that $H$ is covered by the elements of $M$. Let $A \in M$, $A \ne T^D$, and
      suppose that $a \in A^{\times}$ fixes some affine point $z'$; then, as $H$ is a linear group now, $zz'$ is fixed
      pointwise, hence $z$ linewise, so that $a \in T^D$. This implies that $a \in T^D \cap A = Z$, meaning that $a$ is trivial,
      contradiction. It follows that $A$ is a Singer group of the affine plane. Note that doing the same argument on
      the affine lines not incident with $z$, immediately yields the fact that if $A \ne T, T^D$, $A$ is of type (b). 
So each element of $M$ except $T$ and $T^D$
gives a Singer group of $\mP(q)$ under the mapping $\chi$, all not coming from the standard elation group of $\hW(q)$ (and all having the same action on the lines 
incident with $x$).
If $q$ is even, $H$ contains precisely two maximal (elementary) abelian groups of order $q^2$, namely $T$ and $T^D$ \cite{Hira}. So only the standard example arises here.\\

Letting $z$ vary on the line at infinity, we obtain all Singer groups $\overline{K}$ of $\mP(q)$ for which $\mathbb{S} \leq \overline{K}$ and $\overline{K}/\mathbb{S}$
is abelian. So in particular, those examples with $\mathbb{S} \leq \overline{K}$, and $[\overline{K},\overline{K}] \leq \mathbb{S}$. 
The Heisenberg groups are examples of such groups, and so the aimed at classification of \cite{SDWKT2} falls under this description.\\

\subsection{Even characteristic}
Let $q$ be an even prime power, and suppose $q \geq 5$; then by Theorem \ref{Stropp} the automorphisms of $\mP(q)$ are induced by automorphisms of $\hW(q)$ fixing $x$.
Suppose $\overline{K}$ is a Singer group of $\mP(x)$ for which $\overline{K}/(\mathbb{S}\cap \overline{K})$ is abelian.
Then as the latter group is point-transitive on $\Pi(x)$, it is sharply transitive, and whence $\mathbb{S}$ is contained in $\overline{K}$. 
Also, by \cite{DemOst}, we have that there is only one possibility for $\overline{K}/\mathbb{S}$; it is the translation group of $\Pi(x)$. So $\overline{K}$ is unique, as it is 
the translation group of $\hW(q)$ for the point $x$.\\

Denote by $\mathbf{S}^{\mathbb{S}}(\mP(q))$ the set of Singer groups of $\mP(q)$ which are induced by sharply transitive groups of $\Pi(x)$, that is, those {\em containing $\mathbb{S}$}. Then, as in odd characteristic, we also have: 

\begin{theorem}
There is a natural bijection\\

\begin{equation}
\label{lab2}  
\chi:  \mathbf{S}(\Pi(x))  \mapsto  \mathbf{S}^{\mathbb{S}}(\mP(q)) : K \mapsto \chi(K).
\end{equation}
\eop \\
\end{theorem}

Singer groups coming from non-abelian quotients in both characteristics will be handled simultaneously in the next paragraph.\\

\subsection{Non-abelian quotients}

Let $K$ be a Singer group of $\mP(q)$ inducing a sharply transitive group on the points of $\Pi(x)$, that is, let
$K \in \mathbf{S}^{\mathbb{S}}(\mP(q))$. We have classified those $K$ (or $\overline{K}$) for which the quotient $\overline{K}/\mathbb{S}$ is abelian.
In this paragraph, we carry out a classification for the non-abelian case. When $(p,h) = 1$, where $q = p^h$ and $p$ is prime, the classification is complete.
If $(p,h) = p$, an extra assumption will be needed which amounts to demanding that the group be linear.\\

First consider $\overline{K}$; since it is a $p$-group, it fixes some line $L \I x$, whence $\overline{K}/\mathbb{S} = T$ fixes a flag at infinity of $\Pi(x)$, say $(\ell,[\infty])$, where
$[\infty]$ is the line at infinity. We want to consider those $T$ which are contained in the group $H(\ell)$ which is generated by the translation group $A$ ``of'' $[\infty]$ and the translation group $B$
``of'' $\ell$; note that this group is isomorphic to the classical semifield group previously described. (If $T$ is abelian, it is always contained in this group.) We will see in the next paragraph that this assumption is natural, and that in many cases it is not even an {\em extra} assumption.\\

We first note that $H(\ell)/A$ acts sharply transitively on $[\infty] \setminus \{\ell\}$ and that $H(\ell)/B$ acts sharply transitively on the lines of $\Pi(x)$ incident with $\ell$ and different from $[\infty]$. We also note that $Z(H(\ell)) = A \cap B$.

\begin{lemma}
$T \cap B = Z(H(\ell))$.
\end{lemma} 

{\em Proof}.\quad
Follows from the fact that each element of $B$ fixes some line on $\ell$ pointwise, together with $\vert B \cap T\vert = q$ (which follows from $BT = H(\ell)$).\eop \\

The reader notes that $H(\ell)/Z(H(\ell))$ is elementary abelian (so a vector space over $\mathbb{F}_p$).
The following criterion is now immediate.

\begin{theorem}
\label{thminj}
There is a natural injection
$\eta$ from the set $\mathbf{B}(\ell)$ of subgroups  of $H(\ell)/Z(H(\ell))$ of order $q$ which intersect trivially with $B/Z(H(\ell))$, and 
$\mathbf{S}^{\mathbb{S}}(\mP(q))$.
\end{theorem}

{\em Proof}.\quad
The map $\eta$ is defined by the fact that for each $T \in \mathbf{S}(\Pi(x))$ which fixes $\ell$,  we have a short exact sequence

\begin{equation}
\mathbf{1} \mapsto Z(H(\ell)) \mapsto T \mapsto U \mapsto \mathbf{1}, 
\end{equation}

\noindent
where $U \in \mathbf{B}(\ell)$, and vice versa, together with the existence of $\chi$. \eop \\

Note that when $q$ is odd, precisely $q$ elements of $\mathbf{B}(\ell)$ give rise to elements of $\mathbf{S}^{\mathbb{S}}(\mP(q))$ with abelian quotients in $\mathbf{S}(\Pi(x))$ (one of these
elements yielding the translation group of $\Pi(x)$ and so the elation group of $\hW(q)^x$); when $q$ is even, precisely one element does.

The first equation displayed in the next corollary can be found in \cite[Chapter 3]{Hirsch}.

\begin{corollary}
\label{cornumb}
We have that 
\begin{equation}
\vert \mathbf{B}(\ell)\vert = p^{h^2} \leq \vert \mathbf{S}(\mP(q))\vert.
\end{equation}
The number of Singer groups with nonabelian quotients of $\mP(q)$ coming from  $\eta(\mathbf{B}(\ell))$ is
\begin{equation}
p^{h^2} - q^{p\mod{2}}.
\end{equation}
\eop \\
\end{corollary}

Letting $\ell$ vary on $[\infty]$, and noting that  for $\ell \ne \ell'$
precisely one element of $\mathbf{B}(\ell)$ yields a Singer group  coming from some element in $\mathbf{B}(\ell')$, namely the one giving the translation group of $[\infty]$,
we have that  the total number of Singer groups with nonabelian quotients of $\mP(q)$ coming from  $\eta(\bigcup_{\ell}\mathbf{B}(\ell))$ is

\begin{equation}
(q + 1)(p^{h^2} - 1) - {((q + 1)(q - 1))}^{p\mod{2}}.
\end{equation}

\medskip
\subsection{Linear Singer groups, and the case $(p,h) = 1$}

Suppose that $q \geq 5$; then each Singer group of $\mP(q)$ is induced by a subgroup of $\mathrm{Aut}(\hW(q))_x$.
Let $K$ be a {\em linear} Singer group of $\mP(q)$ (that is, induced by some subgroup of $\mathbf{PGSp}_4(q)_x$).

If $q$ is odd, by Theorem \ref{thmoddcent} we have $\mathbb{S} \leq \overline{K}$, so that $\overline{K}/\mathbb{S}$ is sharply transitive and linear on the plane $\Pi(x)$.
Now note that, if $\Omega$ is a Desarguesian projective plane over $\mathbb{F}_q$, and $V$ is any line, then the linear automorphism group $A$ of $\Omega$ which fixes 
$V$ is isomorphic to

\begin{equation}
P \rtimes \PGL_2(q),
\end{equation}
where $P$ is the subgroup of $\Aut(\Omega)$ which fixes $V$ pointwise, and $A/P$ induces the projective general linear group on $\PG(1,q) \cong V$.
So  $\vert A \vert = q^3(q + 1)(q - 1)^2$, and any Sylow $p$-subgroup of $A$ ($q$ being a power of $p$) obviously has the property that 
\begin{itemize}
\item
it fixes some point of $v \I V$;
\item
it has size $q^3$.
\end{itemize}
Vice versa, it is easy to see that for each point $v \I V$, $\mathbf{H}(v) \cong \mathbf{H}_1(q)$ is the unique Sylow $p$-subgroup in $A_v$ (since it is a normal subgroup of $A_v$). \\

It follows readily that the map $\eta$ defines a bijection from $\cup_{\ell}\mathbf{B}(\ell)$ to the linear Singer groups of $\mP(q)$, 
and Corollary \ref{cornumb} gives the precise numerical information.\\

If $(p,h) = 1$, any Singer group of $\mP(q)$ is linear, so that we have a complete classification in odd characteristic.

\begin{theorem}
When $q = p^h \geq 5$ and $p$ is odd, all linear Singer groups are essentially known; in particular, if $(p,h) = 1$, we have a complete classification of Singer groups of 
$\mP(p^h)$.\\
If $q$ is even and $q \geq 8$, all linear Singer groups $K$ of $\mP(q)$ for which $\mathbb{S} \leq \overline{K}$ are known.
\eop \\
\end{theorem}

Note that in all these cases, $\overline{K}/\mathbb{S}$ is naturally embedded in a Heisenberg group $\mathbf{H}_1(q)$, so it has class $\leq 2$.
Whence $K$ has class at most $3$.

We state the prime case in a separate subsection.

\medskip
\subsection{Prime case}

In the next theorem, we ignore the case $q = p = 2$ for obvious reasons.

\begin{theorem}
Let $p$ be an odd prime, and let $K$ be a Singer group of $\mP(p)$. Then $\mathbb{S} \leq \overline{K}$ and $\overline{K}/\mathbb{S}$ is abelian. 
Moreover, there are precisely $p$ different such groups $K$ which fall under two isomorphism classes (exactly one is elementary abelian and the others 
are isomorphic to $\mathbf{H}_1(p)$).
\eop \\
\end{theorem}

The possible actions $\overline{K} \curvearrowright \mS$ were all described earlier in this section.

\bigskip



\section{Central lemma}

As all symplectic GQs naturally embedded in $\mathbf{PG}(3,q)$ are projectively equivalent, we may assume without loss of generality that we are dealing with $\hW(q)$ arising from the symplectic polarity $\rho$ determined by the matrix

\begin{equation} P:=\left(\begin{array}{cccc}  0 & 1 & 0 & 0 \\  -1 & 0 & 0 & 0 \\ 0 & 0 & 0 & 1 \\ 0 & 0 & -1 & 0  \end{array}\right). \end{equation}

We are interested in the behavior of the group of symmetries about a point $x$ of $\hW(q)$ within the stabilizer of that point, $\mathrm{Aut}(\hW(q))_x$.
As the automorphism group of $\hW(q)$ is transitive on the points of $\hW(q)$ we may assume that $x$ is the point $(1,0,0,0)$.
A linear automorphism of $\mathbf{PG}(3,q)$ will fix $(1,0,0,0)$ and stabilize its image $X_1=0$ under $\rho$ if and only if it is induced by a non-singular matrix of the form

\begin{equation} A:=\left(\begin{array}{cccc}  1 & a & b & c \\  0 & d & 0 & 0 \\ 0 & e & f & g \\ 0 & h & i & j  \end{array}\right). \end{equation}

\noindent
Furthermore, every such matrix determines a unique linear automorphism of $\mathbf{PG}(3,q)$ satisfying the above properties.

A linear automorphism of $\mathbf{PG}(3,q)$ determined by a matrix $M$ will induce an automorphism of $\hW(q)$ if and only if 

\begin{equation}
M^TPM=kP,\  \ k\in\mathbb{F}_q.
\end{equation} 

Hence $\mathrm{Aut}(\hW(q))_x$ within $\mathbf{PGL}_4(q)$ is induced by the group 

\begin{equation}G:=\{A\mid A=\left(\begin{array}{cccc}  1 & a & b & c \\  0 & d & 0 & 0 \\ 0 & e & f & g \\ 0 & h & i & j  \end{array}\right), \ A\ \mathrm{nonsingular},\ A^TPA=kP,\ k\neq 0\}.\end{equation}

\noindent
A simple computation also shows that the matrices in $G$ satisfy $d=fj-ig$, and hence a matrix from $G$ has determinant $d^2$.

One easily checks that the group of symmetries about $(1,0,0,0)$ is isomorphic to the subgroup $\mathbb{S}$ of $G$ determined by matrices of the form

\begin{equation} A:=\left(\begin{array}{cccc}  1 & a & 0 & 0 \\  0 & 1 & 0 & 0 \\ 0 & 0 & 1 & 0 \\ 0 & 0 & 0 & 1  \end{array}\right), \ a\in\mathbb{F}_q. \end{equation}

We want to find the largest subgroup $H$ of $G$ such that $\mathbb{S}\leq Z(H)$, that is, $H = C_G(\mathbb{S})$.

One readily checks that in order for a symmetry to commute with a general element of $G$, it is necessary and sufficient that $d=1$. The matrices in $G$ with this property form a subgroup $H$ of $G$.

\medskip
\begin{proposition}
\label{propcent}
\begin{itemize}
\item[{\rm (i)}]
If $q$ is even $H$ is isomorphic to $\mathrm{Aut}(\hW(q))_x\cap\mathbf{PSL}_4(q) = \mathrm{Aut}(\hW(q))_x\cap\mathbf{PGL}_4(q)$.
\item[{\rm (ii)}]
If $q$ is odd, $H$ is a subgroup of index $2$ in $\mathrm{Aut}(\hW(q))_x\cap\mathbf{PSL}_4(q)$.
\end{itemize}
\eop \\
\end{proposition}

Note that for general $q$,
\begin{equation}
\vert \PSL_n(q) \vert = \frac{\vert \PGL_n(q)\vert}{(n,q - 1)}.
\end{equation}

By Proposition \ref{propcent}, any linear Singer group of $\mP(q)$ in even characteristic (seen as an automorphism group of $\hW(q)$) is centralized by the group of symmetries with center $x$.

Now suppose that $K$ is a Singer group of $\mP(q)$, where $q$ is an odd prime power (of the prime $p$). Since the index of $H$ in  $\mathrm{Aut}(\hW(q))_x\cap\mathbf{PGL}_4(q)$
is a power of $2$, the fact that $\PSL_4(q) \unlhd \PGL_4(q)$ implies that $K$ is a subgroup of $H$. (If $K$ would not be in $\PSL_4(q)$, we would have

\begin{equation}
\langle \PSL_4(q),K \rangle = \PSL_4(q)K,
\end{equation}
so that $p$ would divide $[\PGL_4(q) : \PSL_4(q)]$, contradiction. Since $H$ is a normal subgroup of $\PSL_4(q)$ of index $2$, the same reasoning 
forces $K$ to be inside $H$.)\\

We have proved the following theorem.

\begin{theorem}
\label{thmcent}
Let $q$ be any prime power.
Any linear Singer group of $\mP(q)$ is centralized by $\mathbb{S}$.
\eop
\end{theorem}



\bigskip
\section{Isomorphisms}

Since one of our goals is to gain insight in the set of isomorphism classes of Singer groups acting on one and the same quadrangle, and since our setting 
is closely related to central extensions of abelian groups, one might wonder whether a cohomological argument could help us (in at least getting a grip on lower bounds of the number of isomorphism classes). Before passing to concrete 
results on isomorphism classes, we want to take a closer look at this idea.\\

Let $q = p^n$ be a power of the prime $p$, $n \in \mathbb{N}$.
Let $G$ be an elementary abelian group of order $q$, and $A$ a trivial elementary abelian $G$-module of order $q$.
We want to determine the second cohomology group $H^2(G,A)$ by using the Universal Coefficient Theorem (a good reference on this matter is \cite[Chapter 6]{HilSta}), which states that for any natural $i$ we have a short exact sequence

\begin{equation}
0  \longrightarrow \mathrm{Ext}(H_{i - 1}(G,\mathbb{Z}),A) \longrightarrow H^i(G,A) \longrightarrow  \mathrm{Hom}(H_i(G,\mathbb{Z}),A) \longrightarrow 0,
\end{equation}
\noindent
which splits (but not naturally),
where $\mathrm{Ext}(.,.) = \mathrm{Ext}^1_{\mathbb{Z}}(.,.)$. For $i = 2$, we have that $H_1(G,\mathbb{Z}) = G/G' \cong G$, and $H_2(G,\mathbb{Z})$ is the Schur 
multiplier $M(G)$ of $G$. Whence

\begin{equation}
\mathrm{Ext}(H_{1}(G,\mathbb{Z}),A) \cong \mathrm{Ext}(G,A) = \mathrm{Ext}(\oplus_1^{n}\mathbf{C}_p,A) = \bigoplus_1^{n}\mathrm{Ext}(\mathbf{C}_p,A), 
\end{equation}
\noindent
so that
$\mathrm{Ext}(H_{1}(G,\mathbb{Z}),A) \cong \mathbf{C}_p^{n(n - 1)}$. 

Since $G \cong \mathbf{C}_p^{n}$, we have that 
\begin{equation}
M(G) \cong \mathbf{C}_p \oplus \mathbf{C}_p^2 \oplus \cdots \oplus \mathbf{C}_p^{n - 1} = \mathbf{C}_p^{n(n - 1)/2}.
\end{equation}

Let $g_i$ be a generator of the $i$th copy of $\mathbf{C}_p$ in  $\mathbf{C}_p^{n(n - 1)/2}$, $i = 1,\ldots,n(n - 1)/2$. Then a homomorphism $\psi: \mathbf{C}_p^{n(n - 1)/2} \longrightarrow A$
is completely determined by its images of the $g_i$ in $A$, and any mapping 

\begin{equation}
f: \{g_1,\ldots,g_{n(n - 1)/2}\} \longrightarrow A 
\end{equation}
\noindent
yields such a homomorphism $\psi_f$.
It follows that  $\mathrm{Hom}(H_2(G,\mathbb{Z}),A) \cong \mathbf{C}_p^{n^2(n - 1)/2}$, and hence

\begin{equation}
H^2(G,A) \cong \mathbf{C}_p^{n(n - 1)} \oplus   \mathbf{C}_p^{n^2(n - 1)/2}.
\end{equation}

\medskip
Let $q$ be any prime power.
If it is odd, any  linear Singer group of $\mP(q)$ can, by the above, be constructed as a semidirect product 

\begin{equation}
S = T \rtimes \mathbb{S},
\end{equation}
\noindent
where $\mathbb{S} \leq Z(S)$ is the group of symmetries with center $x$. (If $q$ is even, we will assume, besides linearity, that the Singer group contains $\mathbb{S}$.)
Moreover, in the notation of the previous section, $Z(H(\ell)) = B \cap T$ is elementary abelian of order $q$, and  $T/Z(H(\ell))$ also is.\\

Now since we can regard $Z(H(\ell)) = Z$ as a trivial $T/Z$-module, we have that the maximal number of 
such nonisomorphic $T$ is 

\begin{equation}
\vert H^2(T/Z,Z) \vert = p^{n(n - 1)(n + 2)/2}.
\end{equation}

Now consider the map 

\begin{equation}
\varphi: \mathbf{S}(\Pi(x)) \longrightarrow H^2(T/Z,Z): T \rtimes \mathbb{S} \longrightarrow [T].
\end{equation}

(Note that we fix $\ell$ throughout this section without loss of generality.) Now suppose $n > 1$, noting that $\vert H^2(T/Z,Z) \vert = 1$ when $n = 1$.
Let us also not consider abelian groups $T$ (that is, elements of $\varphi^{-1}(\id)$). 
Then the average size of a fiber $\varphi^{-1}([T])$, $[T] \ne \id$,  is at least

\begin{equation}
\frac{p^{n^2} - {q}^{p\mod{2}}}{p^{n(n - 1)(n + 2)/2} - 1} =: F(p,n).
\end{equation}

But $F(p,n) < 1$, so this approach does not lead us to a satisfying result. 

\medskip
\subsection{Odd characteristic}

For now, we suppose that $q$ is odd.
We first prove that, once known that $T \not\cong T'$ in our construction, then $T \rtimes \mathbb{S} \not\cong T' \rtimes \mathbb{S}$ (the 
converse is not necessarily true). 

\begin{lemma}
\label{lemTs}
Let $S$ be a linear Singer group of $\mP(q)$, $q$ odd. Then
$Z(S) = \mathbb{S}$.
\end{lemma}
{\em Proof}.\quad
We already know that $\mathbb{S} \leq Z(S)$ by Theorem \ref{thmoddcent} and Theorem \ref{thmcent}, so suppose $z \in Z(S) \setminus \mathbb{S}$. 
Let $U$ be a line which meets the line $[\ell] \I x$ at infinity (= in $\hW(q)$), which is fixed by $S$; then $\vert S_U\vert = q$.  Note that since $\mathbb{S}$ centralizes
$S_U$, each element of $S_U$ fixes $\proj_Ux = u$ linewise. 
If $z$ moves $U$ to some line not incident with $u$, then each point of $u^{\langle z\rangle}$ is fixed linewise by
$S_U$. This is not possible since $\vert u^{\langle z\rangle}\vert \geq p$, where $q$ is a power of the odd prime $p$, and any line of $\hW(q)$ is antiregular.
So $z$, and then also $\langle z,\mathbb{S}\rangle$, fixes all points of $[\ell]$. But then some nontrivial element of the center must fix some line of $[\ell]^{\perp}$ which is not incident with $x$, 
implying that it fixes all lines of $[\ell]^{\perp}$ (by the action of $S$), again a contradiction. \eop \\

\begin{theorem}
\label{Ts}
Let $q$ be odd. If $T \not\cong T'$, then $T \rtimes \mathbb{S} \not\cong T' \rtimes \mathbb{S}$.
\end{theorem}

{\em Proof}.\quad
By the previous lemma, $S/\mathbb{S} \cong T$ is an invariant of any linear Singer group of $\mP(q)$. \eop \\

\medskip
\subsection{Even characteristic}

Obviously, the lemma in the previous paragraph does not hold in general when $p = 2$. Our approach is much more simple here due to the 
different nature of $\langle A,B\rangle = H(\ell)$. We will come back to this matter at the end of the next section.\\

\bigskip

\section{Partition function and isomorphism classes}

 By Theorem \ref{Ts}, we know that finding nonisomorphic $T$s yields nonisomorphic Singer groups (in the particular setting of that theorem). It seems that classifying these groups could be a very hard problem. We provide a lower bound on the number of isomorphism
classes, but probably the real value is (much) higher than the one we obtain, cf. the discussion below. In even characteristic, things become much easier.

We emphasize that in this section, we will only consider/construct $T$s which live inside the Heisenberg group $H(\ell)$ acting on the associated affine plane $\Pi(s)$.
As such, the obtained Singer  groups (obtained after central extension by the factor $\mathbb{S}$) have at most class $3$ (where class $3$ is the generic case). 
By the nature of Lemma \ref{lemTs}, we have to make a distinction between odd and even characteristic.

\subsection{Odd characteristic}

We will play the following game: we consider a fixed (projective) $(2n - 1)$-space $\eta$ over $\mathbb{F}_p$, $n \in \mathbb{N}_{0,1}$ and 
$p$ an odd prime, and a Desarguesian $(n - 1)$-spread $D = \{ \pi_0,\ldots,\pi_{p^n}\}$ (in fact, the Desarguesian property will not be used, so any 
spread will do). We then want to find out how many different multisets

\begin{equation}
S(\alpha) = \{\vert \alpha \cap \pi_0 \vert,\ldots,\vert \alpha \cap \pi_{p^n}\vert\}
\end{equation}
exist, where $\alpha$ is an $(n - 1)$-space not in $D$. We emphasize that the sets are multisets (taking multiplicities into account), and that we ignore empty intersections.

If $\alpha$ and $\alpha'$ are $(n - 1)$-spaces for which $S(\alpha) \ne S(\alpha')$,
then in the Heisenberg group $H$ corresponding to $\eta$ (i. e., when taking the quotient of $H$ by the center and passing to 
projective space we get $\eta$), to $\alpha$ and $\alpha'$ correspond nonisomorphic groups $T(\alpha)$ and $T(\alpha')$ of order $p^{2n}$; if 
\begin{equation}
[.,.]: H/Z(H) \times H/Z(H) \longrightarrow \mathbb{F}_{p^n} 
\end{equation} 
is the nonsingular bilinear form which comes from the commutator operator, where $H/Z(H)$ is seen as a vector space over $\mathbb{F}_{p^n}$ and $H' = Z(H)$
is identified with $\mathbb{F}_{p^n}$, 
then 
$ \{\alpha \cap \pi_0,\ldots,\alpha \cap \pi_{p^n}\}$ and $ \{\alpha' \cap \pi_0,\ldots,\alpha' \cap \pi_{p^n}\}$ 
describe the commutation structure of $T(\alpha)$ and $T(\alpha')$, so if $S(\alpha)$ and $S(\alpha')$ are different, then
$T(\alpha)$ and $T(\alpha')$ cannot be isomorphic. (The $p^n + 1$ spread elements correspond to the $p^n + 1$ maximal abelian subgroups of order $p^{2n}$ of $H$.)
The function 
\begin{equation}
\zeta(p,n) := \vert \{ S(\alpha) \vert \alpha\  \mbox{is}\ \mbox{an}\ (n - 1)-\mbox{space}\ \mbox{not}\ \mbox{in}\ D \} \vert 
\end{equation}
yields a lower
bound for the number of nonisomorphic $T$s, and hence nonisomorphic Singer groups through our Theorem \ref{Ts}.

We will call the vector $\nu(T(\alpha)) := (\alpha \cap \pi_1,\ldots,\alpha \cap \pi_{p^n})$ the {\em commuting vector} of $T(\alpha)$.\\

\textsc{Connection with the number of directions}.\quad
The reader already notes the difficulty of the problem by making the following consideration: let $\alpha$ be as before, and construct a $2n$-space $\iota$ in which $\eta$ is embedded as a hyperplane. Then from $D$ arises a  (classical)
translation plane in the standard way, and to $\alpha$ correspond $p^n$ disjoint point sets of size $p^n$.
Then $\vert S(\alpha)\vert$ is the number of directions (of the plane) determined by any of these sets. 
But the problem we have to solve is much harder than ``just'' determining the spectrum of directions; we also 
need the intersection information with the spread elements.\\

In any case, we have the following theorem:

\begin{theorem}
The number of nonisomorphic Singer groups of $\mP(q)$, $q = p^n$ odd and $p$ prime, is bounded below by the number of directions of a linear set (over $\mathbb{F}_p$) of size  $q$ in $\mathbf{AG}(2,q)$.
\eop
\end{theorem}

\textsc{Partition function.}\quad
Define $P_{n - 1}$ to be the set of $(n - 1)$-spaces in $\eta$ which are not contained in $D$, and let $P(n)$ be the set of all nontrivial partitions of $n$
(trivial meaning the partition $n = n$).

We suspect that there is an injection

\begin{equation}
\gamma: P(n) \longrightarrow P_{n - 1}: p_n \longrightarrow \gamma(p_n),
\end{equation}
which maps $p_n =  n_1 + n_2 + \cdots  + n_k$ to an $(n - 1)$-space $\gamma(p_n)$ with the property that there exist
spread elements $\pi_{i1},\ldots,\pi_{ik}$ such that $\pi_{ij} \cap \gamma(p_n)$ is an $(n_j - 1)$-space,
and these $k$ intersections generate $\gamma(p_n)$. (In other words, $p_n$ is associated with $n$ independent points such that $\pi_{ij}$ contains 
precisely $n_j$ of these points.)

Note that if this were true, we have found a good lower bound for the possible number of multisets $S(\alpha)$, namely 
the number of nontrivial partitions of $n$.  
For suppose $m$ is a positive integer.
Let $P^+(m)$ be the number of different (unordered) partitions of $m$ including the trivial one; by the Hardy-Ramanujan estimate \cite{HardyRam}, we have the exponential estimate
when $m$ tends to $+\infty$:

\begin{equation}
P^+(m) \longrightarrow \frac{\mathrm{exp}(\pi\sqrt{2m/3})}{4m\sqrt{3}}.
\end{equation}

A generating function for $P^+(m)$ is provided by Euler's reciprocal function:

\begin{equation}
\sum_{v = 0}^{\infty}P^+(v) = \prod_{k = 1}^{\infty}\frac{1}{1 - x^k}.
\end{equation}

\bigskip
\begin{remark}
{\rm
It is easily seen that the cases $n = 2$ and $n = 3$ satisfy our prediction.}
\end{remark}

\subsection{Even characteristic}

For $p = 2$, the problem considered above is a lot simpler; there, we only have to consider two skew $(n - 1)$-spaces $\pi$ and $\pi'$
in a fixed $(2n - 1)$-space $\eta$ over $\mathbb{F}_2$, and consider the number of possible multisets
$\{ \alpha \cap \pi,\alpha \cap \pi' \}$ with $\alpha$ an $(n - 1)$-space different from $\pi$ and $\pi'$. (Since 
the corresponding Heisenberg group has precisely two maximal abelian subgroups of order $p^{2n}$.) 

Note that $\mathbf{P\Gamma L}_{2n}(p)$ acts transitively on the pairs $(\pi,\pi')$, so that any choice is good.

The number we seek is precisely the number of partitions of size $2$ of $n$, that is, 
\begin{equation}
\lceil \frac{n - 1}{2} \rceil,
\end{equation}
which, via the reasonings explained above, is a lower bound for the number of nonisomorphic Singer groups of $\mP(q)$ in even characteristic:

\begin{theorem}
The number of nonisomorphic linear Singer groups of $\mP(q)$ in even characteristic is bounded below by   
\begin{equation}
\lceil \frac{\log_2q - 1}{2} \rceil.
\end{equation}
\end{theorem}

{\em Proof}.\quad
From each of the partitions considered above, a group $T$ arises which lives inside $H(\ell) = \langle A,B \rangle$ (recall that $A$ and $B$ are the only two maximal abelian subgroups of $H(\ell)$), and which contains $A \cap B$. Note that $T \cap A$ and $T\cap B$ are precisely the maximal abelian subgroups of $T$.
Extending by $\mathbb{S}$, we obtain a linear Singer group $S$ of $\mP(q)$
for which $Z(S) \geq \mathbb{S}$ holds. We have to find a variation of Theorem \ref{Ts} which applies to this particular class of Singer groups.
Let $G$ be the translation group of $\hW(q)$ with translation point $x$, and $G^*$ the dual translation group with translation line $L \I x$. 
Note that $S \cap G$ is an abelian group containing $\mathbb{S}$ for which $(S \cap G)/\mathbb{S} = T \cap A$, and that $S \cap G^*$ is an abelian
group containing $\mathbb{S}$ such that $(S \cap G^*)/\mathbb{S} = T \cap B$. It follows easily that $S \cap G$ and $S \cap G^*$ are precisely the maximal abelian 
subgroups of $S$. So once the commuting vector of another group $T' \leq H(\ell)$ such as $T$ is essentially different from that of $T$, the Singer groups 
which arise indeed cannot be isomorphic.
\eop \\

\bigskip
\begin{remark}
{\rm
We notice again that all considerations about  the isomorphism classes of $\mathbf{S}(\Pi(x))$ carry over without any change to general finite semifield planes.\\
}
\end{remark}

\bigskip
\section{Payne-derivatives of TGQs of even order}

A {\em translation generalized quadrangle} (TGQ) is an EGQ with abelian translation group. \\

Every TGQ of order $s$ can be constructed from a {\em generalized oval} in $\mathbf{PG}(3n-1,q)$, that is, a set of $q^n+1$ mutually disjoint $(n-1)$-dimensional subspaces each three of which span $\mathbf{PG}(3n-1,q)$ \cite{TGQ}. Let $\mathscr{O}$ be a generalized oval in $\mathbf{PG}(3n-1,q)$. Then each element $\pi$ of $\mathscr{O}$ is contained in a unique $\mathbf{PG}(2n-1,q)$ that intersects the elements of $\mathscr{O}$ only in $\pi$, the so-called ``tangent space'' of $\mathscr{O}$ at $\pi$. Now suppose that $q$ is even. Then it is well-known that the $q^n+1$ tangent spaces of $\mathscr{O}$ have an $(n-1)$-dimensional space $\nu$ as their mutual intersection, the so-called {\it kernel} or {\em nucleus}. Then 

\begin{equation}
\mathscr{H} := \mathscr{O}\cup\{\nu\}
\end{equation}
 is a {\em generalized hyperoval} in $\mathbf{PG}(3n-1,q)$, that is, a set of $q^n+2$ mutually disjoint $(n-1)$-dimensional spaces, each three of which generates $\mathbf{PG}(3n-1,q)$.  In the case $n=1$ a generalized (hyper)oval simply is a (hyper)oval in $\mathbf{PG}(2,q)$. As of today, the only known generalized (hyper)ovals in $\mathbf{PG}(3n-1,q)$, $n >1$, arise by ``blowing up'' a (hyper)oval in $\mathbf{PG}(2,q^n)$ (that is, by interpreting such a (hyper)oval over $\mathbb{F}_q$).
\medskip

From now on let $q$ be even, $\mathscr{O}$ be a generalized oval in $\mathbf{PG}(3n-1,q)$, and $\mathscr{H}$ the generalized hyperoval arising from it. It is well known that the point $(\infty)$ of the TGQ $\hT(\mathscr{O})$ is a regular point, and hence we can construct the Payne-derivative $\mathscr{P}:=\mathscr{P}(\hT(\mathscr{O}),(\infty))$. This GQ is most easily described as follows. Embed $\mathbf{PG}(3n-1,q)$ as the hyperplane at infinity in $\mathbf{PG}(3n,q)$. 
\begin{itemize}
\item
Then the \textsc{points} of $\mathscr{P}$ are the points of $\mathbf{PG}(3n,q)\setminus\mathbf{PG}(3n-1,q)$, and 
\item
the \textsc{lines} of $\mathscr{P}$ are the $\mathbf{PG}(n,q)$'s that intersect $\mathbf{PG}(3n-1,q)$ exactly in an element of $\mathscr{H}$.
\item
\textsc{Incidence} is the natural one. 
\end{itemize}

We are interested in Singer groups of $\mathscr{P}$. 
\medskip

Clearly the group of all translations of the affine space $\mathbf{PG}(3n,q)\setminus\mathbf{PG}(3n-1,q)$ acts as a Singer group on $\mathscr{P}$. This is an elementary abelian group, and it arises naturally from (and is isomorphic to) the translation group of $\hT(\mathscr{O})$. Conversely every abelian Singer group arises in this way (see \cite{SDWKT1}). 

Of course the main question now is which other groups can act as a Singer group on $\mathscr{P}$. Because the only known generalized ovals are essentially the ovals in $\mathbf{PG}(2,q)$ we will restrict ourselves to the case $n=1$, and discuss the known hyperovals. Also, we will mostly restrict ourselves to linear groups.

For the various hyperovals we will either show that the above mentioned elementary abelian Singer group is the unique Singer group, or provide some non-abelian examples.
We only consider infinite classes.\\

\subsection{The regular hyperoval}

In this case $\mathscr{P}$ is isomorphic to $\mathscr{P}(\hW(q),r)$, where $r$ is any point of $\hW(q)$. Singer groups of this GQ were discussed in the previous section.\\

\subsection{Translation hyperovals}

Let $q=2^h$. Consider the translation hyperoval 

\begin{equation}
\mathscr{H}:=\{(1,t,t^{2^k}) \vert \mathrm{gcd}(k,h)=1 \} \cup \{(0,0,1)\} \cup \{(0,1,0)\}, 
\end{equation}
in $\mathbf{PG}(2,q)$. Note that every translation hyperoval is isomorphic to such a hyperoval. We suppose that $k\neq1, h-1$, so that $\mathscr{H}$ is not a regular hyperoval. In \cite{OKeefe-Penttila1} it was shown that the linear automorphism group of $\mathscr{H}$ is isomorphic to $C_2^h \rtimes C_{q-1}$.  As a Singer group will only contain elements of order a power of $2$ we only need to consider the normal $C_2^h$ part of this automorphism group. This elementary abelian $2$-group of order $q$ acts sharply transitively on the points of $\mathscr{H}\setminus\{(0,0,1),(0,1,0)\}$ and fixes the line $X=0$ pointwise. Its elements are 

\begin{equation}
\left(\begin{array}{ccc} 1 & 0 & 0\\ a & 1 & 0\\ a^{2^k}& 0 & 1  \end{array}\right),\ a \in \mathbb{F}_q.
\end{equation}

Now embed $\mathbf{PG}(2,q)$ as the hyperplane $U=0$ at infinity in $\mathbf{PG}(3,q)$.
It is natural to wonder whether one can ``combine'' this group with the elementary abelian group of order $q^2$ consisting of the elations of $\mathbf{PG}(3,q)$ with axis $U=0$ and center on $X=U=0$ to form a Singer group.  
One easily checks that 

\begin{equation}
G=\langle\left(\begin{array}{cccc} 1 & 0 & 0 & a \\ a & 1 & 0 & b \\ a^{2^k} & 0 & 1 & c \\ 0 & 0 & 0 & 1 \end{array}\right) \vert a,b,c\in\mathbb{F}_q\rangle
\end{equation}
has order $q^3$, acts sharply transitively on the points of $\mathbf{PG}(3,q)\setminus\mathbf{PG}(2,q)$, and stabilizes $\mathscr{H}$ in $\mathbf{PG}(2,q)$. Hence $G$ is a Singer group of the GQ $\mathbf{T}_2^*(\mathscr{H})$. The group $G$ has exponent $4$, a center $Z$ of order $q^2$ (the above mentioned elations), and $G/Z$ is elementary abelian of order $q$.
Note that this construction of course also works for regular hyperovals. In that case this group appears through a different construction in \cite{BambergGiudici}, Lemma 3.12.\\

\subsection{Monomial hyperovals}
Let $\mathscr{H}$ be a monomial hyperoval that is not a translation oval. Then by \cite{OKeefe-Penttila1} it has a linear automorphism group of order $3(q-1)$ or $q-1$. It follows that in this case the elementary abelian translation group is the only linear Singer group. The known examples of such hyperovals are the Segre hyperoval in $\mathbf{PG}(2,2^h)$, $h\geq5$ odd, and the Glynn hyperovals in $\mathbf{PG}(2,2^h)$, $h\geq5$ odd. As in both cases $h$ is odd, it is impossible to have a semilinear Singer group, and hence the elementary abelian Singer group is the unique Singer group in these cases. \\


\subsection{The Payne hyperovals}

These hyperovals exist in all planes $\mathbf{PG}(2,q)$, $q=2^h$, $h\geq5$ odd.  They have a linear automorphism group of order $2$ (see \cite{OKeefe-Penttila1}). As $h$ is odd, every Singer group has to be linear.\\ 

We will now describe a Singer group of $\mathbf{T}_2^*(\mathscr{H})$ distinct from the translation group.
Let $\mathscr{H}$ be the Payne hyperoval  described by $\mathscr{H}=\{(1,t,t^{1/6}+t^{3/6}+t^{5/6}), t\in \mathbb{F}_q\}\cup\{(0,0,1),(0,1,0)\}$. The unique nontrivial linear automorphism of this hyperoval is induced by 

\begin{equation}\gamma:=\left(\begin{array}{ccc} 0 & 1 & 0\\ 1 & 0 & 0\\ 0 & 0 & 1  \end{array}\right).\end{equation}

As before embed $\mathbf{PG}(2,q)$ as the hyperplane $U=0$ at infinity in $\mathbf{PG}(3,q)$. Consider the automorphism 

\begin{equation}g:=\left(\begin{array}{cccc} 0 & 1 & 0 & 1 \\ 1 & 0 & 0 & 0 \\ 0 & 0 & 1 & 0 \\ 0 & 0 & 0 & 1 \end{array}\right).\end{equation}

This automorphism has order $4$ and induces $\gamma$ on $\mathbf{PG}(2,q)$. Next, let $T$ be a subgroup of the translation group of order $2^{3h-1}$ in the following way:

\begin{equation}T=\{h_{(a,b,c)}:=\left(\begin{array}{cccc} 1 & 0 & 0 & a \\ 0 & 1 & 0 & b \\ 0 & 0 & 1 & c \\ 0 & 0 & 0 & 1 \end{array}\right), (a,b,c)\in D\},\end{equation}
where $D$ is an additive subgroup of order $2^{3h-1}$ of $(\mathbb{F}_q,+)^3$ not containing $(1,0,0)$, but containing $(1,1,0)$, and having the property that $(a,b,c)\in D$ implies that $(b,a,c)\in D$. Such $D$ can easily be constructed as follows. Consider $(\mathbb{F}_q,+)^3$ as the vector space $V(3h,2)$, in the obvious way. Then the hyperplane 

\begin{equation}
X_1+\cdots+X_h+X_{h+1}+\cdots X_{2h}=0
\end{equation}
of $V(3h,2)$ defines a subgroup $D$ as required. Now it is easy to check that $g$ normalizes $T$ in $\mathbf{GL}_4(q)$. Hence, as $g^2=h_{(1,1,0)}\in T$, we have that $\left |\langle g\rangle\cap T\right |=2$, and consequently that $\left |\langle g,T\rangle\right |=2^{3h}$. 
Furthermore, since

\begin{equation}h_{(a,b,c)}g=\left(\begin{array}{cccc} 0 & 1 & 0 & 1+a \\ 1 & 0 & 0 & b \\ 0 & 0 & 1 & c \\ 0 & 0 & 0 & 1 \end{array}\right),\end{equation}
we deduce that the point $(0,0,0,1)$ of $\mathbf{PG}(3,q)$ can be mapped onto any other point of $\mathbf{PG}(3,q)\setminus\mathbf{PG}(2,q)$ by an element of $\langle g,T\rangle$ (as $\langle(1,0,0),D\rangle=(\mathbb{F}_q,+)^3$ ). Consequently the group $S:=\langle g,T\rangle$ acts as a Singer group on $\mathbf{T}_2^*(\mathscr{H})$. It is a group of exponent $4$, but its most remarkable property is that it has an intersection with the translation group of maximal possible size $q^3/2$. 

It is important to notice that, as $\gamma$ is a translation of $\mathscr{H}$, this construction can also be applied to the regular hyperoval, and to any other  translation hyperoval for that matter.\\

\subsection{Cherowitzo hyperovals}
These hyperovals exist over all fields of order $2^h$, $h\geq5$ odd. Concerning its automorphism group the following is known. In \cite{OKeefe-Penttila2} it was shown that the linear automorphism group of this hyperoval for $h=5$ is trivial, and in \cite{OKeefe-Thas} it was shown that every linear automorphism that fixes either $(0,0,1)$ or $(0,1,0)$ has to be the identity. No other linear automorphisms of these hyperovals are known, and hence no new Singer groups (that is, distinct from the elementary abelian one) of the related GQ arise for now. \\

\subsection{Subiaco hyperovals}
These hyperovals exist over all fields $2^h$, $h\geq 6$. (For smaller $h$ they fall under classes discussed above.) \\

$\bullet$\quad 
If $h\equiv2\pmod4$ there are two non-isomorphic Subiaco hyperovals in $\mathbf{PG}(2,2^h)$.  The first one has a full (semilinar) automorphism group of order $10h$. The linear automorphism groups contains an elation fixing two points of the hyperoval, and hence a construction analogous to the one described for the Payne-hyperovals can be applied to construct a Singer group of exponent $4$ with an intersection of size $2^{3h-1}$ with the translation group. 
The second Subiaco hyperoval has a full automorphism group of size $5h$. If $h$ is odd it follows that the elementary abelian Singer group is the only one. If $h$ is even there may be semilinear Singer groups arising.  \\

$\bullet$\quad 
If $h\not\equiv2\pmod4$ then there is, up to isomorphism, a unique Subiaco hyperoval. It has a full automorphism group of order $2h$. As above, there is an elation fixing two points of the hyperoval and stabilizing the hyperoval, so that here as well a construction as for the Payne-hyperovals can be carried out to construct a non-abelian Singer group. Again, if $h$ is odd, every Singer group has to be linear, however, for even $h$ there might be semilinear Singer groups as well.\\

\subsection{Adelaide hyperovals}

These hyperovals exist over all fields of order $2^h$, $h\geq6$ even. (For smaller $h$ they fall under classes discussed above.) The linear stabilizer of these hyperovals has order two, whereas the full (semilinear) stabilizer has order $2h$. Again one can construct a Singer group as in the case of the Payne-hyperovals. Here as well there might be semilinear Singer groups arising.\\


\bigskip
\section{Construction of lattices in $\widetilde{\mathbf{C}_2}$-buildings}

Recall that a {\em $\widetilde{\mathbf{C}_2}$-building} is an affine building related to a Coxeter group of type $\widetilde{\mathbf{C}_2}$. It corresponds
to the Coxeter diagram below:\\

\begin{center}
$\widetilde{\mathbf{C}_2}$: \begin{tikzpicture}[style=thick]
\foreach \x in {1,2,3}{
\fill (\x,0) circle (2pt);}
\draw (1,0.035) -- (2,0.035); 
\draw (1,-0.035) -- (2,-0.035); 
\draw (2,0.035) -- (3,0.035); 
\draw (2,-0.035) -- (3,-0.035); 
\end{tikzpicture}
\end{center}

\medskip
An affine building of rank $n \geq 3$ has ``at infinity'' a spherical building of rank $n - 1$. Using this observation, and applying his monumental classification
of spherical buildings of rank at least $3$ \cite{TitsBN},
Tits classified all irreducible affine buildings of rank at least $4$ in \cite{Titsaff}. In the rank $3$ case, such a classification is not possible due to free constructions \cite{Ronan}.
Several other constructions of ``exotic'' rank $3$ buildings are known.

Below, we will describe a recent construction of Essert of (lattices in) $\widetilde{\mathbf{C}_2}$-buildings, taken from \cite{Essert1,Essert2}, which uses finite Singer quadrangles.
One of the advantages of Essert's approach is that his constructions come with nice group actions, and easy presentations.\\

\subsection{Uniform lattices}

Let $G$ be a locally compact topological group. Then there is a left invariant measure 
on G which is unique up to a multiplicative constant, called the {\em Haar measure}. 
This measure induces an invariant measure on every quotient space of $G$. 
A {\em lattice} $\Gamma 
\leq G$ is a discrete subgroup of $G$ such that the measure 
space $\Gamma/G$ has finite measure. 
A {\em uniform lattice} $\Gamma$ is a discrete and cocompact subgroup of $G$, meaning that 
$\Gamma/G$ is  compact. 

Let $X$ be a locally finite piecewise Euclidean complex. Its full automorphism 
group $G = \Aut(X)$ is a topological group, and a neighbourhood basis for the identity is given 
by the stabilizers of compact sets. The vertex stabilizers $G_x$ are then open and compact, 
since 

\begin{equation}
G_x = \lim_{r \longrightarrow \infty} 
{G_x}_{\vert B_r(x)}, 
\end{equation}
\noindent
where $B_r(x)$ is the closed ball around $x$ with radius $r$. Obviously, since $X$ is locally finite, every 
restriction ${G_x}_{\vert B_r(x)}$ is a finite group, so $G_x$ is profinite and hence compact. In particular, 
$G$ is locally compact. We assume additionally that $G$ acts cocompactly on $X$, that is, $X/G$ is compact.

\begin{lemma} 
A subgroup 
$\Gamma \leq G$ is discrete if and only if all its stabilizer subgroups 
of vertices are finite. It is cocompact in $G$ if and only if it acts cocompactly on $X$. 
\end{lemma}

In particular, $\Gamma\leq G$ is a uniform lattice if and only if it acts cocompactly with finite 
stabilizers on $X$.\\

\subsection{Essert's construction} 

Let $G(Y)$ be a developable complex of groups with finite vertex groups 
over a finite ``scwol'', as in \cite{Essert1,Essert2}. Then the geometric realization of the universal cover $\chi$ is a locally 
finite piecewise Euclidean complex. 
Its automorphism group $\Aut(\chi)$ is then a locally compact group with compact open 
stabilizers. The fundamental group $\pi_1(G(Y))$ acts on $\chi$, the stabilizers are the vertex 
groups and the quotient is the geometric realization of the scwol 
$Y$. So $\pi_1(G(Y))$ is a  uniform lattice in $\Aut(\chi)$. 
Hence, to construct lattices in the automorphism groups of locally finite buildings, all 
one has to do is 
construct a finite scwol with finite vertex groups whose universal cover is a 
Euclidean building. \\

For buildings of type 
$\widetilde{\mathbf{C}_2}$, Essert gives two different constructions of panel-regular lattices, 
depending on the types of panels the lattice acts regularly (= sharply transitively) on.

We will describe Essert's construction in a slightly more general setting than was done in \cite{Essert1,Essert2}.\\

Start from a finite generalized quadrangle $\mS$ of order $(s,t)$, $s \ne 1 \ne t$; we write $J = 
\{0,1,\ldots,t\}$ and suppose a Singer group $S$ acts 
sharply transitively on the points of the quadrangle. 
Fix the set of lines $\cL$ through a point $x$ and a 
bijective enumeration function 

\begin{equation}
\lambda : J 
\longrightarrow \cL. 
\end{equation}

We repeat this construction for a second 
quadrangle $\mS'$ of order $(s',t)$ and obtain a Singer group $S'$ and functions $\lambda'$. 

We assume that $\lambda$ and $\lambda'$ are such that for every $j \in J$, we have 
\begin{equation}
S_{\lambda(j)} \cong S'_{\lambda'(j)},
\end{equation}
\noindent
in other words, we have a ``local isomorphism'' between the groups $S$ and $S'$.

In the theorem below, $X \ast Y$, where $X$ and $Y$ are groups, denotes the group generated by $X \cup Y$ modulo all relations
in $X$ and $Y$. 

\medskip
\begin{theorem}[J. Essert \cite{Essert1,Essert2}]
The finitely presented groups 
\begin{equation}
\Gamma_1 = (S \ast S')/\langle [S_{\lambda(j)},S'_{\lambda'(j)}] \vert j \in J \rangle 
\end{equation}
as well as 
\begin{equation}
\Gamma_2 = (S \ast S' \ast \langle c\rangle)/\langle c^{q+2}, c_j \psi_j(x)c_j\psi'_j(x)^{-1} \vert j 
\in J, x \in \mathbb{Z}/q\rangle 
\end{equation}
are uniform panel-regular lattices in a building of type $\widetilde{\mathbf{C}_2}$, which is, respectively, defined as the universal cover $\mB$
of $\Gamma_1$ and $\Gamma_2$. \\
\end{theorem}

If $q > 3$, then the associated buildings are necessarily exotic \cite{Essert1,Essert2}. In any case, the above 
presentations imply very simple descriptions of these buildings. Still, it is not clear when precisely the obtained buildings are actually new.
Of course, as the point- and plane-residues are isomorphic to $\mS$ and $\mS'$, we have some control on the isomorphism class of $\mB$. 
And with respect to the isomorphism class of the constructed lattices, Essert provides the following information on the integral homology.

\begin{theorem}[J. Essert \cite{Essert1,Essert2}]
If $\Gamma_i$ is any of the two lattices as previously constructed, we have 
$H_j(\Gamma_i,\mathbb{Q}) = 0$ for $j  \ne 0$. 
In addition, for the first type of lattices we have
\begin{equation} 
H_2(\Gamma_1,\mathbb{Z}) \cong H_2(S,\mathbb{Z}) \oplus H_2(S',\mathbb{Z}). 
\end{equation}
\end{theorem}

(Note that $H_2(S^{(')},\mathbb{Z})$ is the Schur multiplier of $S^{(')}$.)

\bigskip
\begin{remark}[Property (T)]
{\rm
Note that all cocompact lattices in buildings of type $\widetilde{\mathbf{C}_2}$ have (Kazhdan's) Property (T) if the orders of the GQs as vertex links 
is different from $(2,2)$ \cite{Zuk}.
}
\end{remark}
\medskip

 \subsection{New uniform lattices in $\widetilde{\mathbf{C}_2}$-buildings}

In \cite{Essert1,Essert2}, Essert applies his construction method to the case where $\mS \cong \mS' \cong \mP(q)$, and $S \cong S' \cong \mathbf{H}_1(q)$.
In that case one knows that for all $i$ and $j$, 
\begin{equation}
S_{\lambda(i)} \cong S'_{\lambda'(j)} \cong (\mathbb{F}_q,+).
\end{equation}

Although in general the situation is more subtle,  the constructions also works when $\mS$ and $\mS'$ are allowed to be  nonclassical. We first consider 
the classical Payne-derived quadrangles, admitting not necessarily isomorphic linear Singer groups.

We put 
\begin{equation}
\mS \cong \mP(q) \cong \mS'
\end{equation}
 for some fixed but arbitrary prime power $q = p^n$. We let $S$ and $S'$ be {\em any} 
two linear Singer groups as coordinatized above by elements of the moduli space 
$\{ \mathbf{B}(\ell)\vert \ell \I L\}$. (Properties which hold by this assumption will be freely used below.)

(The only thing we want to suppose is that 
they are not (both) isomorphic to the classical examples, that is, not abelian nor isomorphic to $\mathbf{H}_1(q)$.)\\

 Put $H = H(\ell) = AB$ as before, and consider the Singer group $S =  \mathbb{S} \rtimes T$. (Throughout we suppose that $S$ contains $\mathbb{S}$.)
 
 Put $T \cap A = T_A$ | it is a group that contains $B \cap A$.
 Let $\nu(T)$ be the commuting vector of $T$, and suppose the $p^n + 1$ $(n - 1)$-spaces of the considered spread of $\PG(2n - 1,p)$ are 
 \begin{equation}
 \zeta_0,\zeta_1,\ldots,\zeta_{p^n}.
 \end{equation}
 So $\nu(T) = (\nu(T)_0,\ldots,\nu(T)_{p^n})$, with $\nu(T)_i = \zeta_i \cap \mathbf{P}(T)$ ($\mathbf{P}(T)$ is the projective space corresponding to $T$). We also choose  indexes in such a way that $\zeta_0$ corresponds 
 to $A$ and $\zeta_1$ to $B$.\\
 
 We know that $S$ fixes some line $U \equiv U(\ell)$ incident with $x$ (in the Payne-integrated quadrangle). In fact, it follows from our construction that $S$ acts transitively on the 
lines of $\mS$ incident with $x$ and different from $\ell$. Moreover, since we are working with linear groups, any element of $S$ which fixes some line 
$V \I x$ besides $U$, automatically fixes all lines on $x$. \\

 Now for any line $M \not\sim U$, we have that $S_M$ is an elementary abelian group such that 
 \begin{equation}
 (S_M\mathbb{S}/\mathbb{S})(B \cap A) = T_A, 
 \end{equation}
 and its size is 
 $p^{\mathrm{dim}(\nu(T)_0) + 1}$.

 \begin{proposition}
Let $M$ be any line of $\mS$. 
\begin{itemize}
\item[{\rm (i)}]
If $M \sim U$ in $\hW(q)$, then $S_M \cong (\mathbb{F}_q,+)$.
\item[{\rm (ii)}]
If $M \nI U$ in $\hW(q)$, then $S_M$ is controlled by $\nu(T)$. 
\end{itemize}
\end{proposition}

{\em Proof}.\quad
(i)\quad
Since we are working in $\mathbf{PGSp}_4(q)$, any element of $S_M$ fixes $U$ pointwise, so induces an element of the corresponding translation 
group in $\Pi(x)$. (i) easily follows.\\

(ii) was obtained prior to the statement of the theorem. \eop \\

This has the following implication in the context of Essert's construction. 

 \begin{proposition}[Compatibility]
 A necessary and sufficient condition for Essert's construction to work is ``agreeing on $\nu(.)_0$'' (that is, 
 the necessary bijections on the indexes exist so as to have the local isomorphisms we want).
 \eop \\
 \end{proposition}

 For the sake of convenience, we depict the construction in the following diagram. (In each row, the first two arrows stand for natural projection,
 and the third one for passing to the projective space corresponding to the vector space in question. All upward and downward arrows 
 stand for natural embedding. Also, as before, $\mathbf{P}(R)$ is the projective $\mathbb{F}_p$-space coming from the elementary abelian group $R$.)\\
 

 \begin{equation}
  \hspace*{-1cm}
\begin{array}{ccccccccc}
S = T \rtimes \mathbb{S}& \rightarrow &T &\rightarrow &T/(A \cap B) &\rightarrow &\mathbf{P}(T/(A \cap B)) = \PG(n - 1,p)& &\\
&&&&&&&&\\
 \downharpoonright  &  &               \downharpoonright&    &\downharpoonright&  &\downharpoonright &\searrow&\\
 &&&&&&&&\\
\mathbf{H}_1 \rtimes \mathbb{S}& \rightarrow &H = AB &\rightarrow &H/(A \cap B) &\rightarrow &\mathbf{P}(H/(A \cap B)) = \PG(2n - 1,p)& &S_M, \nu(T)_0\\
&&&&&&&&\\
 \upharpoonright &   &               \upharpoonright&    &\upharpoonright&  &\upharpoonright &\nearrow&\\
 &&&&&&&&\\
\mathbf{H}_1 = A \rtimes \mathbb{S}& \rightarrow &A &\rightarrow &A/(A \cap B) &\rightarrow &\mathbf{P}(A/(A \cap B)) = \zeta_0& &\\
\end{array} 
\end{equation}

\bigskip
\textsc{Odd characteristic}.\quad
When $p$ is odd, we only work with the known $\mS \cong \mS' \cong \mP(q)$. But {\em any} choice of nonisomorphic linear $S$ and $S'$ agreeing on $\nu(.)_0$ leads to new lattices (Essert only considers the  case $S \cong S' \cong \mathbf{H}_1$). \\

\textsc{Even characteristic}.\quad
In the even case, we can also consider nonisomorphic $\mS$ and $\mS'$ (so that the buildings must be different than those constructed by Essert), with respective Singer groups 
$S$ and $S'$ agreeing on $\nu(.)_0$. Then apply Essert.\\


\end{document}